\documentclass{amsart}

\usepackage[ansinew]{inputenc}
\usepackage[french]{babel}
\usepackage{amssymb,amsmath}
\usepackage{graphicx}

\newcommand{\h}{\mathbb H}
\newcommand{\C}{\mathbb C}
\newcommand{\R}{\mathbb R}
\newcommand{\V}{\mathcal V}
\newcommand{\N}{\mathcal N}
\newcommand{\ksimple}[2]{\Lambda^#1_s#2}
\newcommand{\bs}{\backslash}
\newcommand{\kvect}[2]{#1_1\wedge\dots\wedge{#1}_{#2}}

\newcommand{\mymod}{\mbox{ mod }}
\newcommand{\myspan}[1]{<\!{#1}\!>}
\newcommand{\si}{S_\phi}
\newcommand{\leg}{{\mathcal L}}
\newcommand{\fibrunit}{S_\phi TM}

\newcounter{compte}
\newcounter{dessin}

\theoremstyle{definition}

\newtheorem{definition}[compte]{D\'efinition}

\theoremstyle{plain}
\newtheorem{proposition}[compte]{Proposition} 
\newtheorem{lemme}[compte]{Lemme}
\newtheorem{corollaire}[compte]{Corollaire}
\newtheorem{theorem}[compte]{Théorème}

\author{G. Berck}
\address{Gautier Berck, Universit\'e Catholique de Louvain,
Institut de Math\'ematiques Pures et Appliqu\'ees,
2, Chemin du Cyclotron, B--1348,
Louvain la Neuve}
\email{g.berck@math.ucl.ac.be}

\title[Minimalité des sous-variétés totalement géodésiques]{Minimalité des sous-variétés
 totalement géodésiques en géométrie finslerienne}
\begin{document}
\begin{abstract}
Nous montrons que les sous-variétés totalement géodésiques d'une variété de Finsler 
sont minimales pour le volume de Holmes-Thompson. Pour la mesure de Hausdorff, des contre-exemples 
à ce résultat sont connus (voir \cite{AB}).
\end{abstract}
\maketitle
\section{Introduction}
L'étude des notions d'aires et de volume sur les espaces normés et les variétés de Finsler, 
initiée par Busemann, a récemment acquis un regain d'intérêt, bénéficiant notamment des
apports de géométrie convexe, différentielle et intégrale (voir par exemple \cite{A0},\cite{A2},
\cite{BI},\cite{ga}, \cite{sc} et \cite{AB}). Contrairement 
au cas euclidien, plusieurs dé\-fi\-ni\-tions naturelles d'aire et de volume coexistent dont deux 
se distinguent par leurs propriétés intéressantes. D'une part, les espaces normés et variétés 
de Finsler étant métriques, la mesure de Hausdorff est un choix naturel. D'autre part, le 
volume de Holmes-Thompson, introduit dans \cite{ht} pour des considérations purement géométriques, 
apparaît comme la notion adéquate pour généraliser aux espaces normés les formules 
classiques de géométrie intégrale (voir par exemple \cite{A2} et \cite{sc}). Une autre propriété 
renforce 
l'intérêt porté sur celui-ci : le volume de Holmes-Thompson d'une variété de Finsler est, à une constante près, 
le volume de Liouville de son fibré cotangent unitaire. Ceci ouvre la voie aux techniques symplectiques 
pour l'étude des propriétés géométriques de ce volume.

Equipé de notions d'aire, de volume, on est naturellement amené à étudier les sous-variétés 
minimales et, disposant aussi d'une norme, à comparer les problèmes variationnels uni et 
multi-dimensionnels. Des résultats intéressants pour le volume de Holmes-Thompson ont déjà été 
obtenus dans cette direction. Par exemple, dans \cite{A2} \'Alvarez et Fernandes montrent que, si la métrique 
est projective (c'est à dire si les géodésiques sont des droites), alors les 
$k$-densités de Holmes-Thompson sont des densités de Crofton, ce qui a pour conséquence que les $k$-plans 
affins sont minimaux (voir \cite{AF}). Par des techniques différentes, Burago et Ivanov prouvent 
dans \cite{BI} que les 2-disques affins d'un espace normés minimisent le volume de Holmes-Thompson 
parmi tous les disques immergés de même bords. Aussi, le  'filling theorem' d'Ivanov \cite{I1} a 
pour conséquence que les sous-variétés de dimension 2 totalement géodésiques d'une variété de Finsler 
sont minimales pour 
cette aire (voir \cite{AB}), tandis qu'\'Alvarez et Berck prouvent dans \cite{AB} qu'il en est de même 
pour les hypersurfaces totalement géodésiques.

Dans ce papier, nous prouverons le théorème suivant qui s'inscrit naturellement dans ce contexte:

\noindent{\bf Théorème} {\it Toute sous-variété totalement géodésique d'une variété de Finsler 
est minimale pour le volume de Holmes-Thompson.}

En géométrie riemannienne, ce résultat est obtenu facilement. On montre que la courbure 
moyenne d'une sous-variété est la trace de la seconde forme fondamentale qui elle-même 
s'annule pour les sous-variétés totalement géodésiques. 

Ce n'est par contre pas un résultat général en géométrie de Finsler, sous-entendu pour une 'bonne' 
notion de volume, mais est propre à la densité de Holmes-Thomp\-son. En effet, dans \cite{AB}, les auteurs 
présentent une famille de métriques de Finsler projectives pour lesquelles les plans affins, évidemment 
totalement géo\-dé\-siques, ne sont pas tous minimaux pour la mesure de Hausdorff. De ce fait, aucun 
des résultats cités précédemment n'est valide pour cette mesure.

Il n'y a pas, en géométrie de Finsler, de seconde forme fondamentale adaptée au problème variationnel. 
Néanmoins, pour toute densité $\phi$ et toute sous-variété $N$, on construira en chaque 
point de celle-ci un {\em covecteur de courbure moyenne} $h$ avec la propriété que, pour toute 
variation lisse de $N$, $\alpha(s)$, telle que $\partial_s\alpha(0)$ est nul en dehors d'un
compact de $N$,
$$
\frac{d}{ds}\left(\int_{\alpha(s)}\phi\right)_{\vert s=0}=\int_Nh(\partial_s\alpha)\phi.
$$
Il en résulte qu'une sous-variété est minimale si et seulement si ce covecteur est nul en tout 
point.

La construction de ce covecteur passe par la généralisation aux densités des {\em applications de 
Legendre} et {\em formes de Hilbert} bien connues en géométrie de Finsler. La principale différence 
avec ce cas classique tient dans le fait que, pour une densité donnée, ces applications et formes 
ne sont pas uniques. Il apparaîtra toutefois que le covecteur de courbure moyenne est bien 
défini.

Dans \cite{14}, Busemann et al. construisent, en chaque point de la sphère unité de la densité, une forme 
linéaire particulière que nous appellerons {\em forme de Busemann}. Celle-ci leur permit de montrer que la 
densité de Holmes-Thompson est localement convexe; de notre côté, elle servira à définir une application 
de Legendre et, par suite, une forme de Hilbert particulièrement bien adaptées à cette densité. La proposition \ref{tool}, 
point fondamental dans la preuve du résultat principal quoique quelque peu technique, montrera que cette forme de 
Hilbert est étroitement liée à la forme de Hilbert associée à la norme. Puisque le covecteur de courbure moyenne 
est directement construit à partir de la forme de Hilbert, 
cela nous donnera en retour un lien direct entre les problèmes variationnels uni et multi-dimensionnels.

Ce papier s'articule suivant quatre sections. Nous commencerons par un rapide 
survol de la notion de densité et plus particulièrement de celle de Holmes-Thompson 
dans les espaces normés, suivi de la construction de la forme de Busemann et de la preuve {\it actualisée} (et 
largement inspirée de \cite[pp.18-20]{BI}) de la locale 
convexité de cette densité. Le lemme de structure terminera cette partie sur 
une propriété du fibré unitaire attaché à une densité.

Nous introduisons dans la seconde partie les outils fondamentaux que sont les 
applications de Legendre et formes de Hilbert associées. 

La section suivante sera consacrée au problème variationnel et à la notion de
covecteur de courbure moyenne. Nous effectuerons ensuite le calcul explicite des 
covecteurs de courbure moyenne pour les courbes et surfaces de $\R^3$ et montrons 
qu'on récupère dans ces cas la courbure des courbes et la courbure moyenne des 
surfaces. Nous verrons aussi le lien entre ce covecteur de courbure et la courbure 
des courbes d'une surface de Finsler telle que l'a définie Cartan dans \cite{C1}.

Nous donnerons dans la dernière section la preuve de la proposition liant forme de Hilbert attachée à la norme 
et forme de Hilbert attachée à la densité de Holmes-Thompson, ensuite celle du théorème principal.
Enfin, nous présenterons deux classes 
intéressantes de métriques de Finsler possédant des sous-variétés totalement géodésiques.

L'auteur aimerait particulièrement remercier J.C. \'Alvarez pour les nombreuses discussions fructueuses, 
ses conseils judicieux et pour lui avoir présenté ce problème de sous-variétés minimales; ainsi que P. Lambrechts 
pour ses multiples relectures attentives et ses remarques constructives.

\section{Densités}
\subsection{Densités}
Considérons un espace vectoriel réel $V$ de dimension $n$ et soit $k\leq n$. 
L'ensemble des $k$-vecteurs simples de $V$, $\ksimple{k}{V}$, est un cône dans $\Lambda^kV$ 
appelé {\em cône de Grassmann}. Tout $k$-vecteur simple non nul $a=\kvect{v}{k}$ porte 
un unique $k$-plan vectoriel orienté de $V$, $span\{v_1,\dots,v_k\}$, que nous noterons 
$\myspan{a}$. Bien entendu, si $t$ est un réel positif, $\myspan{ta}=\myspan{a}$. En 
conséquence, le cône de Grassmann est un cône sur la grassmannienne des $k$-plans vectoriels 
orientés de $V$. Nous noterons $\myspan{a}^*$ le $k$-espace vectoriel orienté dual de $\myspan{a}$.
\begin{definition}
Une $k$-densité $\phi$ sur un espace vectoriel réel $V$ de dimension $n$ ($1\leq k\leq n$) 
est une application $\phi\colon\ksimple{k}{V}\to[0,\infty)$ lisse sur $\ksimple{k}{V}\bs\{0\}$,
positive ($\phi(a)=0\Leftrightarrow a=0$) et homogène de degré 1 ($\phi(ta)=|t|\phi(a)$).
\end{definition}
Remarque. Classiquement, la positivité n'est pas exigée et notre définition correspond aux {\em $k$-densités de 
volume lisses} de \cite{A3}. Nous perturbons quelque peu la terminologie par simplicité, puisque nous 
ne considérerons que ce type de densités. La notion de densité est à mettre en rapport avec celle 
d'intégrand paramétrique définie par Federer, voir \cite[p. 515]{FE}. La principale différence, mais non l'unique, 
est qu'un intégrand paramétrique est défini sur tout l'espace vectoriel $\Lambda^kV$. Aussi, une densité est-elle 
la restriction au cône de Grassmann d'un intégrand paramétrique lisse, positif et pair.

La {\em sphère unité} $S_\phi$ est l'ensemble des $k$-vecteurs simples pour lesquels
$\phi=1$. Notons que celle-ci est difféomorphe à la grassmannienne des $k$-plans orientés 
de $V$.

Géométriquement, les densités généralisent les notions d'aires et de volumes des 
espaces euclidiens. Le réel positif $\phi(\kvect{v}{k})$ s'interprète comme le 
volume du parallélotope engendré par les vecteurs $v_1,\dots,v_k$. De plus, tout 
comme les formes différentielles, on peut intégrer les densités sur les sous-variétés 
(non nécessairement orientées) et donc parler de $k$-aire d'une sous-variété compacte 
de dimension $k$ pour une certaine $k$-densité.

La notion de densité s'étend naturellement aux variétés. On reprend, sur chaque
espace tangent, la construction faite sur un espace vectoriel. On considère alors le fibré 
en cônes $\ksimple{k}{TM}\to M$ de fibre $\ksimple{k}{T_xM}$ au dessus de $x$.
\begin{definition}
Soit $M$ une variété de dimension $n$. Une $k$-densité sur $M$ est une application
$\phi\colon \ksimple{k}{TM}\to[0,\infty)$ lisse hors de la section nulle et telle que
pour tout $x\in M$, l'application $\phi_x:=\phi(x,.)\colon \ksimple{k}{T_xM}\to[0,\infty)$ 
est une $k$-densité sur l'espace vectoriel $T_xM$.
\end{definition}
Pour simplifier les notations, on notera $\pi$ la projection $\pi\colon \Lambda^kTM\to M$ ainsi que
ses différentes restrictions : au fibré en cônes $\pi\colon \ksimple{k}{TM}\to M$ et au fibré
unitaire $\pi\colon \fibrunit\to M$ où $\fibrunit=\{(x,a)\in\ksimple{k}{TM}|\phi(x,a)=1\}$.

Nous nous intéresserons plus particulièrement aux densités sur les espaces de Minkowski : les
espaces normés dont la sphère unité, lisse, a toutes ses courbures principales strictement
positives pour une métrique euclidienne auxiliaire quelconque.
Nous noterons par la suite $B\subset V$ et $S=\partial B$ les boules et
sphères unités de $(V,\Vert\cdot\Vert)$. De même, $B^*$ et $S^*$ désigneront les boules et
sphères unités duales de $(V^*,\Vert\cdot\Vert^*)$.

Classiquement, les espaces euclidiens possèdent, pour toute dimension $k$, une densité naturelle :
celle pour laquelle le $k$-volume de tout cube unitaire de dimension $k$ vaut 1.
La situation se complexifie dans le cas des espaces normés non-euclidiens. Pour ceux-ci, diverses
considérations géométriques ont motivé la construction de différentes densités.
Parmi celles-ci, on trouve la densité de Busemann, donnant la mesure de Haussdorff
des sous-variétés, et la densité de Holmes-Thompson, liée au volume de Liouville. Dans la suite, nous
focaliserons notre attention sur cette dernière.

Cette densité est apparue dans différents contextes, avec des définitions différentes mais
équivalentes (voir \cite{A3} pour les différentes définitions). Nous utiliserons ici la définition suivante
qui remonte à Busemann et qu'il appelait alors {\it fonction de projection} \cite{14}.
\begin{definition}
Soit $(V,\Vert\cdot\Vert)$ un espace de Minkowski.
Soit $a\in\ksimple{k}{V}$ et considérons la projection $P_a\colon V^*\to\myspan{a}^*$ duale
de l'injection $\myspan{a}\hookrightarrow V$. La {\em densité de Holmes-Thompson} est définie par
$$
\phi(a)=\frac{1}{\epsilon_k}\int_{P_a(B^*)}a,
$$
où $\epsilon_k$ est le volume euclidien de la boule unité euclidienne de dimension $k$.
\end{definition}

\subsection{Forme de Busemann}
Dans \cite{hl}, Harvey et Lawson introduisent la notion de calibration d'une densité par une forme
différentielle fermée pour étudier les sous-variétés minimales. Nous utiliserons ici une définition 
plus restrictive, équivalente à la notion de convexité (locale) d'une densité (voir \cite{14}). 
\begin{definition}
Une $k$-forme $\xi\in\Lambda^kV^*$ {\em calibre} une $k$-densité $\phi$ en
$a\in\ksimple{k}{V}$ si $\xi(a)=\phi(a)$ et $\xi\leq\phi$ sur le cône de Grassmann.
On dira qu'elle calibre {\em localement} la densité en $a$ si la seconde condition n'est
satisfaite que dans un voisinage de $a$.
\end{definition}
Dans \cite{14}, Busemann et al. montrent que les $k$-densités ($1\leq k\leq n-1$) de Holmes-Thompson
d'un espace de Minkowski de dimension $n$ sont localement calibrées en chaque point de
la sphère unité par une $k$-forme distinguée. Cette forme jouera un grand rôle dans la suite.
Nous reprenons ici leur construction. Nous donnons ensuite, en langage moderne, la preuve de la
calibration locale de la densité de Holmes-Thompson par cette forme (voir \cite[p.24]{14}, \cite[pp.18-20]{BI} 
pour la construction, \cite[p.34]{14} pour la preuve).

Pour $a\in\ksimple{k}{V}$ donné, considérons la projection $P_a\colon V^*\to\myspan{a}^*$ et
soit $\Sigma_a$ le lieu des points singuliers de la restriction $P_a\colon S^*\to\myspan{a}^*$
(voir figure \ref{figu1}). Puisque $S^*$ est quadratiquement convexe,
$\Sigma_a$ est une sous-variété de dimension $k$ de $S^*$ et la restriction
$P_a\colon \Sigma_a\to\partial P_a(B^*)$ est un difféomorphisme. La forme de Busemann se définit comme suit :
\begin{definition}\label{bus}
Soit $V$ un espace de Minkowski et $a\in\ksimple{k}{V}$. La {\em forme de Busemann} en $a$ est définie par
$$
\beta_a(b)=\frac{1}{k\epsilon_k}\int_{\Sigma_a}i_\chi b
$$
où $i_\chi b$ est la contraction de $b$ avec le champ d'Euler $\chi(v)=v$.
\end{definition}
\begin{proposition}[Busemann et al.\cite{14}]
La forme de Busemann calibre localement la densité de Holmes-Thompson.
\end{proposition}
\begin{proof}Considérons $\tilde\chi$ le champ d'Euler sur $\myspan{a}^*$. Un
calcul en coordonnées donne $d(i_{\tilde\chi}a)=ka$. Cette égalité et le théorème de
Stokes appliqués à la définition de la densité de Holmes-Thompson donnent
$$
\phi(a)=\frac{1}{k\epsilon_k}\int_{\partial P_a(B^*)}i_{\tilde\chi}a.
$$
On vérifie ensuite facilement que $P_a^*(a_{|\myspan{a}^*})=a_{|V^*}$ et que
$DP_a(\chi)=\tilde\chi$. En outre, puisque $S^*$ est quadratiquement convexe, la restriction
$P_a\colon\Sigma_a\to\partial P_A(B^*)$ est un difféomorphisme. On obtient donc
$\phi(a)=\beta_a(a)$.

Prenons $b\in\ksimple{k}{V}$ et considérons la projection $P_b\colon V^*\to\myspan{b}^*$. Si $b$ est
suffisamment proche de $a$, alors $P_b(\Sigma_a)$ est une sous-variété de
$\myspan{b}^*$ difféomorphe à une sphère de dimension $k-1$. 
Alors $\beta_a(b)$ est le volume, mesuré avec $b$, du sous-ensemble $A\subset\myspan{b}^*$
de bord $\partial A=P_b(\Sigma_a)$ :
\begin{eqnarray*}
\beta_a(b)&=&\frac{1}{k\epsilon_k}\int_{\Sigma_a}i_\chi b\\
&=&\frac{1}{\epsilon_k}\int_Ab
\end{eqnarray*}
Par la convexité de $B^*$,
$P_b(\Sigma_a)\subset P_b(B^*)$. Donc,
$\beta_a(b)\leq\phi(b)$, (voir figure \ref{figu2}).
\end{proof}

Remarque. Par contre, la densité de Holmes-Thompson n'est en général pas calibrée. En effet,
pour des $k$-vecteurs $c$ éloignés de $a$, il peut arriver que $P_c(\Sigma_a)$ ne soit
pas un plongement de $\Sigma_a$ dans $\myspan{c}^*$. Dans ce cas, les volumes doivent
être comptés avec multiplicités et
on perd en général l'inégalité. Pour des exemples explicites de telles situations où la forme
de Busemann est de plus la seule forme qui calibre localement la densité, voir
\cite{BI}, \cite{14}, (voir figure \ref{figu3}).

\begin{minipage}{4.5cm}
\includegraphics{fig1.1}
\refstepcounter{dessin}\label{figu1}
\centerline{\small {\bf Figure \arabic{dessin}.}}
\end{minipage}
\begin{minipage}[t]{3.5cm}
\includegraphics{fig3.1}
\refstepcounter{dessin}\label{figu2}
\centerline{\small {\bf Figure \arabic{dessin}.}}
La densité de Holmes-Thompson est loca\-le\-ment calibrée ...

\end{minipage}
\hfill
\begin{minipage}[t]{3.5cm}
\includegraphics{fig2.1}
\refstepcounter{dessin}\label{figu3}
\centerline{\small {\bf Figure \arabic{dessin}.}}
... mais en général pas globalement calibrée.

\end{minipage}
\subsection{Lemme de structure}
Pour terminer cette section, nous présentons le lemme de structure. Celui-ci, et sa version duale, nous
seront utiles pour les propositions \ref{hilbert} et \ref{tool}.

Considérons tout d'abord un fibré différentiable $\pi\colon E\to M$.
Nous appellerons {\em fibré vertical} le fibré $\V\to E$ où $\V=\ker D\pi\subset TE$ et
{\em fibré normal} le fibré $\N=TE/\V\to E$.
On notera
$\Gamma(\V)$ et $\Gamma(\N)$ les espaces de sections de ces fibrés. Utilisant l'intégrabilité
du fibré $\V$, on montre aisément que l'application suivante est bien définie et
qu'il s'agit d'une dérivée covariante (voir \cite{Bo} p. 32).
\begin{definition}
La dérivée covariante de Bott est l'application
\begin{eqnarray*}
&\nabla\colon&\Gamma(\V)\times\Gamma(\N)\to\Gamma(\N)\\
&&(X,Y)\mapsto\nabla_XY=[X,\tilde{Y}]\mymod\V
\end{eqnarray*}
où $\tilde{Y}$ est un champ de vecteurs sur $E$ tel que $\tilde{Y}\mymod\V=Y$.
\end{definition}
Comme toute dérivée covariante, celle-ci s'étend à n'importe quel produit tensoriel de $\N$.
Notons
$L_X$ la dérivée de Lie d'un champ de tenseurs par rapport au champ de vecteurs $X$, on a les
dérivées covariantes suivantes :
\begin{eqnarray*}
&&\nabla\colon \Gamma(\V)\times\Gamma(\Lambda^k\N)\to\Gamma(\Lambda^k\N), \nabla_XY=L_X\tilde Y\mymod \V\\
&&\nabla\colon \Gamma(\V)\times\Gamma(\Lambda^k\N^*)\to\Gamma(\Lambda^k\N^*), \nabla_XY=L_X\tilde Y\mymod \V^*
\end{eqnarray*}
En particulier, considérons le fibré $\pi\colon E=\Lambda^kTM\to M$.
Nous noterons, comme précédemment, $\V=ker D\pi$ le fibré vertical et $\N=(T\Lambda^kTM)/\V$ le fibré normal.
La $k^{\grave{e}me}$ puissance extérieure du fibré normal possède une section canonique
$\sigma\colon \Lambda^kTM\to\Lambda^k\N, \sigma(x,a)=b\mymod \V$ avec $D\pi b=a$.
Ceci permet de définir pour tout $(x,a)\in\Lambda^kTM$ une application linéaire
\begin{eqnarray*}
H_{(x,a)}&\colon &\V_{(x,a)}\to\Lambda^kT_xM\\
&&v\mapsto D\pi_{(x,a)}\nabla_v\sigma
\end{eqnarray*}
Puisque le fibré $E\to M$ est vectoriel, l'espace vertical en un point s'identifie naturellement à
la fibre passant par ce point. Nous noterons $I_{(x,a)}\colon \V_{(x,a)}\to\Lambda^kT_xM$ cet
isomorphisme canonique.
\begin{lemme}[Lemme de structure] En tout point $(x,a)\in\Lambda^kTM$, l'isomorphisme
canonique $I_{(x,a)}$ et l'application linéaire $H_{(x,a)}$ co\"\i ncident.
\end{lemme}
\begin{proof}
Soit $(x,a)\in E$, $v$ un vecteur vertical en $(x,a)$ et $\gamma\colon \R\to\pi^{-1}(\pi(x,a))$
une courbe verticale telle que $\gamma(0)=(x,a)$ et $\dot{\gamma}(0)=v$.
Alors, un calcul élémentaire en coordonnées montre que pour toute section $Y\in\Gamma(\N)$,
$$D\pi_{(x,a)}(\nabla_vY)=\frac{d}{dt}D\pi(Y(\gamma(t)))_{\vert t=0}.$$ Le lemme de structure s'en
déduit immédiatement en remarquant que, pour la section canonique $\sigma$, $D\pi(\sigma(x,a))=a$.
\end{proof}
Notons comme conséquence immédiate que si on se restreint au fibré unitaire d'une densité $\phi$,
l'image de l'application linéaire $H_{(x,a)}$ est alors le tangent de la sphère unité
$T_a\si\subset\Lambda^kT_xM$.

On peut faire une construction similaire dans le cas du fibré dual. Nous ne l'utiliserons, dans le
théorème \ref{tool}, que dans le cas du fibré cotangent. Nous ne présenterons donc que ce cas.

Notons $I_{(x,p)}\colon \V_{(x,p)}\to T^*_xM$ l'isomorphisme canonique. Il existe sur $T^*M$ une 1-forme
canonique, $\alpha\in\Omega^1(T^*M), \alpha_{(x,p)}=p\circ D\pi$. Celle-ci s'annule sur le fibré vertical,
on peut donc l'identifier à une section du dual du fibré normal: $\underline{\alpha}\in\Gamma(\N^*)$.
En outre, en utilisant la formule de Cartan, on trouve :
$\nabla_X\underline{\alpha}=(di_X\alpha+i_Xd\alpha)\mymod \V^*$.
Comme $X$ est vertical, on a $\nabla_X\underline{\alpha}=i_Xd\alpha\ \mymod \V^*$.

Comme dans le cas précédent, on construit alors l'application linéaire :
\begin{eqnarray*}
H^*_{(x,p)}&\colon &\V_{(x,p)}\to T^*_xM\\
&&v\mapsto D\pi_{(x,p)}i_vd\alpha
\end{eqnarray*}
On obtient enfin le lemme de structure dual :
\begin{lemme}[Lemme de structure dual]
En tout point $(x,p)\in T^*M$, l'isomorphisme
canonique $I_{(x,p)}$ et l'application linéaire $H^*_{(x,p)}$ co\"\i ncident.
\end{lemme}

\section{Application de Legendre et forme de Hilbert}
\subsection{Application de Legendre}
L'application de Legendre est un outil important du calcul des variations. Dans le cas d'un espace normé $(V,\phi)$, 
de norme lisse, elle est définie par 
$$
\leg\colon V\to V^*,v\mapsto\frac{1}{2}(d\phi^2)_v.
$$
Il ne s'agit pas en général d'un difféomorphisme de $V\bs\{0\}$ vers $V^*\bs\{0\}$. C'est toutefois le 
cas lorsque la forme bilinéaire symétrique définie par sa différentielle en tout vecteur $v$ non nul 
\begin{eqnarray*}
D_v\leg&\colon&T_vV\times T_vV\to\R\\
&&(x,y)\mapsto D_v\leg(x)(y)
\end{eqnarray*}
est définie positive. On dit alors qu'on à affaire à un problème régulier du calcul des variations. Du point 
de vue de la géométrie de Finsler, c'est aussi la définition intrinsèque, et usuelle, d'un espace de Minkowski.

Remarquons que l'application de Legendre est l'unique application $\leg\colon V\to V^*$ homogène de degré 1 
telle que $\leg(v)(v)=\phi^2(v)$ et telle que $\leg(v)$ s'annule sur le tangent $T_{v}S_\phi$ pour tout 
vecteur unitaire $v$. Nous 
utiliserons cette caractérisation pour étendre cette notion aux densités quelconques.
\begin{definition}
Une application de Legendre pour une $k$-densité $\phi$ sur $V$
est une application
$\leg\colon \ksimple{k}{V}\to\Lambda^kV^*$ lisse sur $\ksimple{k}{V}\bs\{0\}$
et telle que
\begin{enumerate}
\item $\mathcal L$ est homogène de degré 1,
\item $\leg(a)(a)=\phi(a)^2$,
\item $\leg(a)$ s'annule sur $T_a\si$.
\end{enumerate}
\end{definition}
Dans le cas $k=n-1$, le cône de Grassmann est l'espace $\Lambda^kV$ tout entier et
$\si$ est une hypersurface de celui-ci.
Il n'existe alors qu'une seule application de Legendre donnée par
$\leg(a)=\frac{1}{2}(d\phi^2)_{a}$. Dans le cas général, il en existe par contre une infinité et il convient alors 
d'en déterminer une bien adaptée à la densité. Pour Holmes-Thompson, c'est la forme de Busemann
qui va nous fournir une application de Legendre adéquate :
\begin{proposition}\label{theta} Une application de Legendre pour la densité de Holmes-Thompson $\phi$ sur un
espace de Minkowski est donnée par
\begin{eqnarray*}
\leg&\colon &\ksimple{k}{V}\to\Lambda^kV^*\\
&&a\mapsto\phi(a)\beta_a
\end{eqnarray*}
où $\beta_a$ est la forme de Busemann définie en \ref{bus}.
\end{proposition}
\begin{proof} Cette application est lisse car la norme duale d'un espace de Min\-kow\-ski est lisse. On vérifie aussi 
de suite qu'elle est homogène. Ensuite, on sait que $\beta_a(a)=\phi(a)$, donc $\leg(a)(a)=\phi(a)^2$. 
Il reste alors à montrer que $\beta_a$ s'annule sur le tangent de la sphère unité. 
On sait que $\beta_a$ calibre localement la densité $\phi$ en
$a$. Donc, $a\in\si$ est un minimum local de la fonction
$$
\phi-\beta_a\colon \ksimple{k}{V}\to\R.
$$
Puisque $\phi$ est constante sur sa sphère unité $\si$, on en déduit que
$$
\forall b\in T_a\si,\ \beta_a(b)=\partial_b(\phi-\beta_a)=0.
$$
\end{proof}
Remarquons que cette construction se généralise à toute densité localement calibrée.

On généralise naturellement aux densités sur les variétés :
\begin{definition}
Une application de Legendre pour une $k$-densité $\phi$
sur une variété $M$ est une application
$
\leg\colon \ksimple{k}{TM}\to\Lambda^kT^*M
$
lisse hors de la section nulle et telle que
\begin{enumerate}
\item $\mathcal L$ se projette sur l'application identité sur $M$ : $\pi^*\circ\leg=\pi$.
\item en tout point $x$ de $M$, l'application $\leg_x\colon \ksimple{k}{T_xM}\to\Lambda^kT_x^*M$
définie par $\leg(x,a)=(x,\leg_x(a))$ est une application
de Legendre sur l'espace vectoriel $T_xM$.
\end{enumerate}
\end{definition}
\subsection{Forme de Hilbert}
A l'aide de l'application de Legendre, nous généralisons la notion classique de forme 
de Hilbert au cas des $k$-densités. Nous proposons deux définitions équivalentes. La première, 
généralisant la définition classique, facilitera dans la section suivante la preuve de 
l'invariance du covecteur de courbure moyenne vis à vis d'un choix d'application de Legendre; la 
seconde sera utile pour prouver le théorème \ref{tool}.
\begin{definition}
Soit $\leg\colon \ksimple{k}{TM}\to\Lambda^k{T^*M}$ une application de Legendre pour une
$k$-densité $\phi$.
La {\em forme de Hilbert} pour cette application est
une $k$-forme sur le fibré unitaire, $\omega_k\in\Omega^k(\fibrunit)$, définie par
$
(\omega_k)_{(x,a)}=\leg_x(a)\circ D\pi.
$
\end{definition}
Il existe sur l'espace dual $\Lambda^kT^*M$ une $k$-forme canonique $\alpha_k$ définie
par $(\alpha_k)_{(x,\xi)}=\xi\circ D\pi^*$. Ceci donne une définition alternative pour la forme de Hilbert.
\begin{definition}[Alternative]
Soit $\leg\colon \ksimple{k}{TM}\to\Lambda^k{T^*M}$ une application de Legendre. La $k$-forme
de Hilbert associée à cette application de Legendre est définie par
$
\omega_k=(\leg^*\alpha_k)_{|\fibrunit}.
$
\end{definition}
On vérifie immédiatement que ces deux définitions sont équivalentes. 

La proposition suivante est une conséquence du lemme de structure. Rappelons qu'au fibré $\pi\colon\fibrunit\to M$ 
sont associés les fibrés vectoriels vertical $\V=\ker D\pi\to\fibrunit$ et normal $\N=(T\fibrunit)/\V\to \fibrunit$.
\begin{proposition}\label{hilbert}
Soit $\phi$ une $k$-densité et $\omega_k$ la forme de Hilbert associée à une application de Legendre.
Soient $u,v\in \V_{(x,a)}$ deux vecteurs tangents verticaux et
$b\in\Lambda^kT_{(x,a)}\fibrunit$ tel que $D\pi_{(x,a)}b=a$. Alors,
(1) $i_u\omega_k=0$, (2) $i_{u\wedge v}d\omega_k=0$ et (3) $d\omega_k(u\wedge b)=0$.
\end{proposition}
\begin{proof} Le premier point est une conséquence immédiate de la définition de forme de Hilbert. 
Pour le second, il suffit de remarquer que la distribution des espaces tangents à la fibre est intégrable et qu'elle est
contenue dans l'annulateur de $\omega_k$. 
Enfin, le dernier point est une conséquence du lemme de structure.
Le premier point de la proposition assure que la forme de Hilbert s'identifie canoniquement
à une section $\underline{\omega}_k\in\Gamma(\Lambda^kN^*)$. De plus, si $\sigma\in\Gamma(\Lambda^kN)$
représente la section canonique, alors par définition de la forme de Hilbert, $\underline{\omega}_k(\sigma)\equiv1$. 
Dès lors,
\begin{eqnarray*}
0&=&\partial_u(\underline{\omega}_k(\sigma))\\
&=&(\nabla_u\underline{\omega}_k)(\sigma)+\underline{\omega}_k(\nabla_u\sigma)
\end{eqnarray*}
D'après le lemme de structure, $D\pi(\nabla_u\sigma)\in T_a\si$, donc 
$\underline{\omega}_k(\nabla_u\sigma)=0$. Alors $(\nabla_u\underline{\omega}_k)(\sigma)=0$.
Or, $b\mymod \V=\sigma(x,a)$ et $\nabla_u\underline{\omega}_k=i_ud\omega_k\mymod \V^*$. On en
déduit que $i_ud\omega_k(b)=0$.
\end{proof}
\section{Formule de première variation}
Le {\it relevé tangent} d'une sous-variété orientée $N^k\subset M$ est la sous-variété du fibré unitaire
$$
\tilde{N}:=\{(x,a)\in\fibrunit|\myspan{a}=T_xN\}\subset\fibrunit.
$$
Bien entendu, la projection
$\pi\colon \tilde{N}\to N$ est un difféomorphisme. Alors par définition de la forme de Hilbert,
\begin{eqnarray*}
\int_N\phi=\int_{\tilde{N}}\omega_k,
\end{eqnarray*}
quelque soit le choix de la forme. Ceci nous permet d'aborder le problème variationnel
dans le contexte classique des formes différentielles.
\begin{proposition}
Soit $\phi$ une $k$-densité sur une variété $M$ et $\omega_k$ une forme de Hilbert pour cette densité.
La sous-variété $N\subset M$ de dimension $k$ est point critique du problème variationnel si et seulement si
$i_b d\omega_k=0$ quelque soit le $k$-vecteur tangent $b\in\ksimple{k}{T\tilde{N}}$.
\end{proposition}
\begin{proof}
Soit $f\colon N\times[0,1]\to M$ une variation de $N$ à support compact. Elle induit une
variation des relevés tangents : $\tilde{f}\colon \tilde{N}\times[0,1]\to\fibrunit$. Soit
$X=\partial_t$ le champ de vecteurs correspondant sur $\fibrunit$.
Alors
\begin{eqnarray*}
\frac{d}{dt}_{|t=0}\int_{N_t}\phi&=&\frac{d}{dt}_{|t=0}\int_{\tilde{N}_t}\omega_k\\
&=&\int_{\tilde{N_0}}di_X\omega_k+i_Xd\omega_k\\
&=&\int_{\tilde{N_0}}i_Xd\omega_k
\end{eqnarray*}
par Stokes puisque la variation est à support compact. Puisque que la variété $N$ est
point critique si et seulement si ceci est nul pour toute variation, on obtient
le résultat.
\end{proof}
Il apparaît clairement au vu de la preuve précédente que 
la condition de minimalité est indépendante d'un choix de forme de Hilbert.

Plus généralement, la différentielle extérieure de la forme de Hilbert permet de
construire, en tout point d'une sous-variété fixée, un covecteur de courbure moyenne
qui mesure la première variation du volume lorsqu'on déforme la sous-variété.
Pour une sous-variété orientée $N\subset M$,
il est naturel, d'après ce qui précède, de considérer le covecteur
$
i_b d\omega_k
$
où $b\in\ksimple{k}{T\tilde{N}}$ est un $k$-vecteur tangent tel que $\omega_k(b)=1$, 
c'est à dire tel que $D\pi_{(x,a)}b=a$. On vérifie comme précédemment que
ce covecteur ne dépend pas d'un choix de forme de Hilbert.
\begin{proposition}
Soit $N\subset M$ une $k$-sous-variété et $\omega_k,\omega_k^\prime$ deux $k$-formes de Hilbert pour
une densité. Alors,
$i_b d\omega_k=i_b d\omega_k^\prime$ quelque soit le $k$-vecteur
tangent $b\in\ksimple{k}{T\tilde{N}}$.
\end{proposition}
\begin{proof}
Quelles que soient les formes de Hilbert, on a toujours 
$$
\int_N\phi=\int_{\tilde{N}_0}\omega_k=\int_{\tilde{N}_0}\omega_k'.
$$
Alors, pour une variation quelconque, en utilisant les mêmes notations que dans la 
preuve précédente, 
$$
\int_{\tilde{N}_0}i_Xd\omega_k=\int_{\tilde{N}_0}i_Xd\omega_k'.
$$
Puisque cette égalité reste vérifiée pour toute variation, le résultat s'en suit.
\end{proof}

D'après le point (3) de la proposition \ref{hilbert}, ce covecteur s'annule sur l'espace tangent à la fibre de
$\fibrunit$. Il se projette donc naturellement sur un covecteur défini sur $M$.
\begin{definition}
Soit $\phi$ une $k$-densité sur une variété $M$ et $N\subset M$ une $k$-sous-variété orientée de
$M$. Le covecteur de courbure moyenne de $N$ en $x\in N$ est défini par
$$
h_x(u)=d\omega(\tilde{u}\wedge b)
$$
où $\tilde{u}\in T_{(x,a)}\fibrunit$ tel que $D\pi(\tilde{u})=u$ et 
$b\in\ksimple{k}{T_{(x,a)}\tilde{N}}$ tel que $D\pi b=a$.
\end{definition}
\begin{corollaire}
Une sous-variété est minimale si et seulement si son covecteur de courbure moyenne est nul en tout point. 
\end{corollaire}

Nous terminons cette section par quelques illustrations classiques.

{\bf Courbes de $\R^3$} Considérons $\R^3$ avec la métrique standard. Nous noterons 
$(x,y,z,u,v,w)$ les coordonnées d'un point du fibré tangent $T\R^3\cong\R^3\times\R^3$ où 
$(x,y,z)$ sont les coordonnées du point et $(u,v,w)$ celles du vecteur basé en ce point.

On vérifie sans peine que la forme de Hilbert, $\omega_1=(\frac{1}{2}d\Vert\cdot\Vert^2)\circ D\pi$, 
est donnée par $\omega_1=<V\vert D\pi(\cdot)>$ au point $(X,V)$. En coordonnées, il s'agit de 
la restriction au fibré unitaire de la forme $ud\!x+vd\!y+wd\!z$.

Soit $c(t)$ une courbe lisse paramétrisée par longueur d'arc. Son relevé tangent dans le fibré unitaire 
est simplement $\tilde{c}=(c,\dot{c})$, tandis que l'unique vecteur tangent à ce relevé qui 
se projette sur $\dot{c}$ est $U=(\dot{c},k)$ où $k=\ddot{c}$ est le vecteur normal multiplié par 
la courbure.

On calcule alors que $d\omega_1(\cdot\wedge U)$ est la restriction au fibré unitaire de 
la forme $\dot{c}_1d\!u+\dot{c}_2d\!v+\dot{c}_3d\!w-k_1d\!x-k_2d\!y-k_3d\!z$. D'après le point (3) de la 
proposition \ref{hilbert}, sur le fibré unitaire, $\dot{c}_1d\!u+\dot{c}_2d\!v+\dot{c}_3d\!w=0$. On 
obtient donc que le covecteur de courbure moyenne est, comme attendu, $h=<-k\vert\cdot>$.

{\bf Surfaces de $\R^3$} Considérons $\R^3$ avec l'aire euclidienne. Comme d'habitude, nous identifions 
$\Lambda^2T\R^3$ avec $T\R^3\cong\R^3\times\R^3$ en identifiant le 2-vecteur $v\wedge w$ avec le 
produit vectoriel $v\times w$. Nous prendrons cette fois $(x_1,x_2,x_3,v_1,v_2,v_3)$ comme 
coordonnées.

On vérifie de nouveau sans peine que la 2-forme de Hilbert $\omega_2$ est la restriction au 
fibré unitaire de la forme $v_1d\!x_2\wedge d\!x_3+v_2d\!x_3\wedge d\!x_1+v_3d\!x_1\wedge d\!x_2$.

Soit $f(r,s)$ une surface régulière. Son relevé tangent est $\tilde{f}=(f,n)$ où 
$n$ est le vecteur normal unitaire. Le covecteur de courbure moyenne est obtenu par la restriction au 
fibré unitaire de 
$$
d\omega_2(\cdot\wedge\frac{\partial_r\tilde{f}\wedge\partial_s\tilde{f}}{\Vert\partial_rf\wedge\partial_sf\Vert}), 
$$
c'est à dire, en utilisant l'expression $n=\frac{\partial_rf\times\partial_sf}{\Vert\partial_rf\times\partial_sf\Vert}$,
$$
\sum_{i=1}^3n_id\!v_i+\frac{(\partial_rn\times\partial_sf+\partial_rf\times\partial_sn)_i}
{\Vert\partial_rf\times\partial_sf\Vert}dx_i.
$$
De nouveau, le point (3) de la proposition 
\ref{hilbert} assure que, sur le fibré unitaire, $\sum_in_id\!v_i=0$. De plus, 
$\partial_rn=\alpha\partial_rf+\beta\partial_sf$ et $\partial_sn=\gamma\partial_rf+\delta\partial_sf$ 
pour certaines fonctions $\alpha,\beta,\gamma,\delta$. Alors un simple calcul donne 
$$
\frac{\partial_rn\times\partial_sf+\partial_rf\times\partial_sn}
{\Vert\partial_rf\times\partial_sf\Vert}=(\alpha+\delta)n=-Mn
$$
où $M$ est l'opposé de la trace de la différentielle de
l'application de Gauss, c'est à dire la courbure moyenne de $f$. On en déduit que le covecteur de
courbure moyenne est, comme attendu, $h=-M<n\vert\cdot>$.

\subsection{Métrique de Finsler sur une surface}
Dans \cite{C1}, Cartan utilise sa méthode d'équivalence pour étudier les invariants différentiels 
d'une surface munie d'une 1-densité. La forme de Hilbert apparaît naturellement dans sa démarche et
lui permet notamment de définir la notion de courbure d'une courbe. Nous esquissons rapidemment,
dans le cas d'une métrique de Finsler, ses résultats et montrons le lien entre courbure
et covecteur de courbure moyenne comme défini précédemment.

Considérons une surface de Finsler $(M,\phi)$ orientée et rappelons que, par définition, la
différentielle de l'application de Legendre définit pour chaque vecteur unitaire $v\in T_xM$
un produit scalaire $<\cdot|\cdot>_v$ sur $T_xM$.
Il existe donc sur le fibré unitaire deux 1-formes canoniques $\omega_1$ et $\omega_2$ définies
par, pour $X\in T_{(x,v)}\fibrunit$,
$$
D\pi(X)=\omega_1(X)v+\omega_2(X)w
$$
où $(v,w)$ est la base orientée de $T_xM$ orthonormée pour le produit scalaire $<\cdot|\cdot>_v$.

Notons que la forme $\omega_1$ n'est autre que la forme de Hilbert et que les relevés tangents
des courbes satisfont l'équation $\omega_2=0$.

Calculant les différentielles extérieures de ces formes, Cartan montre l'existence d'une
$3^{\mbox{\tiny ème}}$ forme canonique, $\omega_3$, linéairement indépendante des deux autres.
Il prouve aussi que ces formes et leurs différentielles sont reliées par les {\em équations de
structure de Cartan}:
$$
\left\{
\begin{array}{l}
d\omega_1=-\omega_2\wedge\omega_3\\
d\omega_2=\omega_1\wedge\omega_3-I\omega_2\wedge\omega_3\\
d\omega_3=-K\omega_1\wedge\omega_2-J\omega_2\wedge\omega_3
\end{array}
\right.
$$
où $I,\ J$ et $K$ sont des fonctions sur le fibré unitaire caractéristiques de la métrique de Finsler. Il montre
également que les extrémales sont ces courbes dont les relevés tangents sont solutions du système
différentiel $\omega_2=\omega_3=0$, ce qui revient à dire que le covecteur de courbure moyenne est nul
d'après la première équation de structure.

Contrairement aux deux autres, la $3^{\mbox{\tiny ème}}$ forme canonique ne s'annule pas sur l'espace tangent à
la fibre du fibré unitaire. En conséquence, cette forme appliquée sur un vecteur tangent du relevé tangent
d'une courbe donne une mesure de la différence entre le 2-jet de cette courbe et celui de la géodésique de
même vecteur tangent. Notant cela, Cartan définit naturellement la courbure d'une courbe comme le rapport
$k=\frac{\omega_3(V)}{\omega_1(V)}$ où $V$ est une vecteur tangent au relevé tangent de cette courbe. De
la première équation de structure, on déduit que le covecteur de courbure moyenne d'une courbe est
$h=-k\omega_2$.
\section{Minimalité des sous-variétés totalement géodésiques}
Nous avons déterminé, à la proposition \ref{theta}, une application de Legendre particulière pour la densité 
de Holmes-Thompson. Celle-ci définit en retour une $k$-forme de Hilbert que nous noterons 
$\Theta_k$. Nous allons montrer, dans la proposition \ref{tool}, que cette forme et la 1-forme de Hilbert 
associée à la norme sont étroitement liées. Or, comme nous l'avons vu dans la section précédente, la 
minimalité des sous-variétés s'exprime par une condition sur la différentielle extérieure de la 
forme de Hilbert. Dès lors, cette proposition établit un pont entre les problèmes variationnels uni et 
multi-dimensionnels. Nous l'utiliserons comme élément clef dans la preuve du théorème principal.

Les 2 formes $\omega_1\in\Omega(STM)$ et $\Theta_k\in\Omega^k(\fibrunit)$ sont définies sur des variétés 
différentes. Pour les comparer, nous allons les ramener sur la variété $M$.

Considérons à cet effet une section locale
$
\sigma\colon M\to\fibrunit.
$
Celle-ci définit une distribution de $k$-plans tangents orientés sur $M$ : $\{\myspan{\sigma(x)}|x\in M\}$.
Prenant les vecteurs unitaires de ces $k$-plans tangents, 
nous définissons un sous-fibré du fibré unitaire : 
$
\pi_\sigma\colon U_\sigma=\{(x,v)\in STM|\ v\in\myspan{\sigma(x)}\}\to M.
$

Nous allons montrer que le pull-back de $\Theta_k$ par la section $\sigma$ est, à une constante près, le push-forward
par $\pi_\sigma$ du produit extérieur de $\omega_1$ et d'une puissance de $d\omega_1$. Plus précisément,
\begin{proposition}\label{tool}
Soit $\sigma\colon M\to\fibrunit$ une section locale. Alors, les formes de Hilbert $\Theta_k$
pour la $k$-densité de Holmes-Thompson et $\omega_1$ pour la norme sont reliées par l'équation :
$$
\sigma^*\Theta_k=\frac{(-1)^{\frac{k(k+1)}{2}}}{k!\epsilon_k}\pi_{\sigma*}(\omega_1\wedge (d\omega_1)^{\ k-1}).
$$
\end{proposition}
Pour prouver ce théorème, nous utiliserons le lemme suivant (cfr prop 2.3, 2.4 de \cite{A1}).
\begin{lemme}
Soit $(V,\Vert.\Vert)$ un espace de Minkowski, de sphère unité $S$.
Soit $W\subset V$ un sous-espace vectoriel et $\Sigma_W$ le lieu des points singuliers de la projection
$P\colon S^*\subset V^*\to W^*$. Alors,
$
\Sigma_W=\leg^1(W\cap S)
$
où $\leg^1$ est l'application de Legendre associée à la norme.
\end{lemme}
\begin{proof}[Preuve de la proposition \ref{tool}]
La fibre $\pi^{-1}_\sigma(x)=\myspan{\sigma(x)}\cap\ ST_xM$ est une sphère de dimension $k-1$ dans
$\myspan{\sigma(x)}$. Puisque $\myspan{\sigma(x)}$ est naturellement orienté, la fibre $\pi^{-1}_\sigma(x)$
l'est aussi.

Par le lemme précédent, $\leg^1(\myspan{\sigma(x)}\cap\ ST_xM)=\Sigma_{\myspan{\sigma(x)}}\subset ST_x^*M$. De plus,
par définition, $\omega_1=\leg^{1*}\alpha_1$.

Soit $\Sigma:=\bigcup_x\Sigma_{\myspan{\sigma(x)}}$ et $\pi_\Sigma\colon \Sigma\to M$ la restriction à $\Sigma$ de la
projection $\pi\colon T^*M\to M$. Alors,
$$
\pi_{\sigma*}(\omega_1\wedge d\omega_1^{k-1})=\pi_{\Sigma*}(\alpha_1\wedge d\alpha_1^{\ k-1}).
$$
Par définition du push-forward,
$$
\pi_{\Sigma*}(\alpha_1\wedge d\alpha_1^{\ k-1})_x(b)=\int_{\pi_\Sigma^{-1}(x)}i_{\tilde b}(\alpha_1\wedge d\alpha_1^{\ k-1}).
$$
Soit $(x,p)\in \pi_\Sigma^{-1}(x)$. Si $(u_1,\dots,u_{k-1})$ est une base de $(ker D\pi_\Sigma)_{(x,p)}$, alors, 
d'après la proposition \ref{hilbert},
$$
i_{\tilde b}(\alpha_1\wedge d\alpha_1^{\ k-1})(\bigwedge_{j=1}^{k-1}u_j)=
(-1)^{\frac{k(k+1)}{2}}(k-1)!(\alpha_1\wedge(\bigwedge_{j=1}^{k-1}i_{u_j}d\alpha_1)(\tilde b).
$$
Utilisant le lemme de structure dual, on vérifie que ceci est encore égal à
$$
(-1)^{\frac{k(k+1)}{2}}(k-1)!\ b\cdot(p\wedge(\bigwedge_{j=1}^{k-1}I(u_j))).
$$
De plus, $(I(u_1),\dots,I(u_{k-1}))$ est une base de $T_p\Sigma_{\myspan{\sigma(x)}}\subset ST^*_xM$. On obtient alors
$$
\pi_{\Sigma*}(\alpha_1\wedge d\alpha_1^{\ k-1})(b)=(-1)^{\frac{k(k+1)}{2}}(k-1)!\int_{\Sigma_{\myspan{\sigma(x)}}}i_\chi b.
$$
De ce fait,
$$
\frac{(-1)^{\frac{k(k+1)}{2}}}{k!\epsilon_k}\pi_{\sigma*}(\omega_1\wedge d\omega_1^{\ k-1})=\beta_{\sigma(x)}\in\Lambda^kT^*_xM.
$$
Enfin, puisque $\Theta_{\sigma(x)}=\leg(\sigma(x))\circ D\pi=\beta_{\sigma(x)}\circ D\pi$, on
obtient $(\sigma^*\Theta)_x=\beta_{\sigma(x)}$.
\end{proof}

Nous pouvons à présent montrer le résultat principal:
\begin{theorem}
Toute sous-variété totalement géodésique d'une variété de Finsler est minimale pour 
la densité de Holmes-Thompson.
\end{theorem}
\begin{proof}
Soit $x$ un point d'une sous-variété de Finsler $M$ et soit $N\subset M$ la sous-variété (locale) 
formée de toutes les géodésiques passant par $x$ et tangentes à un $k$-plan fixé de $T_xM$. On va 
montrer que la courbure moyenne de $N$ est nulle en $x$.

Prenons $\sigma\colon M\to\fibrunit$ une section locale telle que $\sigma(N)\subset\tilde{N}$.
Alors, si $\sigma(x)=(x,a)$, le covecteur de courbure moyenne en $x$ se calcule comme
suit :
$$
h_x=(\sigma^*d\Theta_k)_x(\cdot\wedge a),
$$
où $\Theta_k$ est la forme de Hilbert particulière associée à la transfomation de Legendre 
définie par la forme de Busemann. Utilisant 
la proposition \ref{tool} et le fait que la différentielle extérieure commute avec le pull-back et 
le push-forward, on obtient 
$$
h_x(u)=\frac{(-1)^{\frac{k(k+1)}{2}}}{k!\epsilon_k}\pi_{\sigma*}(d\omega_1)^k(u\wedge a).
$$
Prenons $v\in ST_xN$ tel que $(x,v)\in S_\sigma$. Alors, $v\wedge a=0$ et donc 
$a=v\wedge v_1\wedge\dots\wedge v_{k-1}$ pour certains $v_i$.

Pour effectuer le push-forward, on est amené à calculer 
$$
(d\omega_1)^{k}_{(x,v)}(\tilde{u}\wedge\tilde{v}\wedge\kvect{\tilde{v}}{k-1}\wedge\kvect{u}{k-1})
$$
où $(u_1,\dots,u_{k-1})$ est une base de l'espace tangent à la fibre $T_{(x,v)}\pi_\sigma^{-1}(x)$ et 
$\tilde{v},\ \tilde{v_i},\ \tilde{u}$ des relevés de $v,\ v_i,\ u$ dans $T_{(x,v)}S_\sigma$. Par antisymétrie, 
c'est indépendant du choix des relevés.

Considérons l'unique géodésique $\gamma$ passant par $x$ et tangente à $v$. Par définition de la variété $N$, 
c'est une courbe de $N$ et donc son relévé tangent $\tilde{\gamma}$ est une courbe de $S_\sigma$. Comme 
relevé de $v$, on peut donc choisir $\dot{\tilde{\gamma}}$. Or $\gamma$ étant une géodésique, elle est minimale 
et donc 
$$
d\omega_1(\dot{\tilde{\gamma}},.)=0.
$$
Puisqu'on peut faire ce raisonnement en tout point $(x,v)$ de la 
fibre $\pi_\sigma^{-1}(x)$, on obtient 
$$
h_x=0.
$$
\end{proof}

Nous terminons en présentant deux classes intéressantes de métriques de Finsler admettant des sous-variétés 
totalement géodésiques, l'une classique, l'autre beaucoup moins.

\subsection{Métriques projectives} 
Dans son étude du $4^{\mbox{\tiny ème}}$ problème de Hilbert, Busemann a construit une large
classe de métriques de Finsler sur $\R^n$ dont les géodésiques sont les droites. Son idée fut
de considérer une mesure positive sur l'ensemble des hyperplans affins et de définir la longueur
d'un segment comme la mesure de l'ensemble des hyperplans qui l'intersectent. Puisque tout hyperplan
coupant un segment coupe aussi toute courbe joignant les extrémités de celui-ci, la longueur du segment
sera nécessairement inférieure à la longueur de la courbe. Plus précisément,
\begin{theorem}[Busemann]
Si $\Phi$ est une mesure lisse positive sur l'espace $\mathcal H$ des hyperplans affins de $\R^n$, alors il
existe une métrique de Finsler $\phi$ sur $\R^n$ telle que pour chaque courbe lisse $\gamma:[a,b]\to\R^n$,
$$
\int_{\xi\in\mathcal H}\#(\xi\cap\gamma)\Phi=\int_\gamma\phi
$$
où $\#(\xi\cap\gamma)$ représente le nombre d'intersection de l'hyperplan $\xi$ avec la courbe $\gamma$. De plus,
les géodésiques de cette métrique de Finsler sont les droites affines.
\end{theorem}
Les espaces affins sont évidemment totalement géodésiques pour de telles métriques et donc minimaux pour la
densité de Holmes-Thompson.

\subsection{Métriques sur $\C P^n$ et $\h P^n$}
Dans leur article sur les submersions en géo\-mé\-trie de Finsler \cite{A1}, \'Alvarez et Dur\'an construisent des métriques
de Finsler sur les espaces projectifs complexes et quaternioniques dont les géodésiques sont des cercles
et pour lesquelles les $\C P^1$ et $\h P^1$, et donc tous les sous-espaces projectifs, sont totalement géodésiques. 
Le grand intérêt de cet exemple par rapport au précédent est que ces métriques n'ont pas toutes les mêmes 
géodésiques. 

Nous présentons leur construction sur les espaces projectifs complexes. Il suffira de remplacer le mot 
'complexe' par 'quaternionique' dans ce qui suit pour obtenir une autre classe de métriques.

\begin{definition}
Une submersion $\pi:M\to N$ entre deux variétés de Finsler est dite iso\-métrique si en tout point de $M$, 
l'image de la boule unité par la différentielle de la projection est la boule unité au point image.
\end{definition}
\begin{theorem}[\'Alvarez, Dur\'an] Soit $\varphi$ une métrique de Finsler sur la sphère 
$S^{2n+1}\subset\C^{n+1}$ invariante sous l'action de $S^1$ et dont les géodésiques sont les grands 
cercles, alors il existe une unique métrique de Finsler $\psi$ sur $\C P^n$ telle que la fibration 
de Hopf, $\pi:(S^{2n+1},\varphi)\to(\C P^n,\psi)$, soit une submersion isométrique. De plus, les géodésiques de 
cette métrique sont toutes des cercles et les $\C P^1$ sont totalement géodésiques.
\end{theorem}

Nous donnerons les grands points de la preuve de ce théorème. 

A l'aide de la construction de Busemann (sur la sphère), les auteurs construisent une classe de métriques 
de Finsler sur la sphère unité de $\C^{n+1}$, invariante par l'action de $S^1$, et dont les géodésiques sont 
les grands cercles. Fixons-en une et remarquons que, du fait de l'invariance sous l'action de $S^1$, les 
boules unités basées en des points d'une même fibre ont toutes la même image par la différentielle de 
la projection. Ceci définit alors une métrique de Finsler sur $\C P^n$, et par définition, la submersion entre ces 
variétés de Finsler est isométrique.

Notons que, comme conséquence immédiate de la définition de submersion iso\-métrique, la différentielle de 
la projection décroît les normes. Forts de ce constat, les auteurs prouvent que les géodésiques de la 
variété image sont les projections de {\em certaines} géodésiques de la variété initiale qu'ils appellent 
{\em géodésiques horizontales}.

Les géodésiques de $\C P^n$ sont donc les projections de certains grands cercles de $S^{2n+1}$. Or tout grand 
cercle étant dans un 2-plan complexe, sa projection est dans un $\C P^1$. Puisqu'ils sont soit disjoints, soit 
transverses, les $\C P^1$ sont donc nécessairement totalement géodésiques. 

Enfin, les auteurs 
prouvent que les géodésiques de ces métriques, nécessairement toutes fermées, sont des cercles.

\end{document}